\documentclass{jdcs-new}

 \currentvolume{18}
 \currentnumber{1}
 \received{July 28 2010}
 \revised{May 18 2011}
\allowdisplaybreaks

\usepackage{cite}
\usepackage{graphicx}

\usepackage[all]{xy}

\theoremstyle{plain}
\newtheorem{theorem}{Theorem}

\newtheorem{proposition}[theorem]{Proposition}

 \theoremstyle{definition}
\newtheorem{definition}[theorem]{Definition}
 \theoremstyle{remark}

\begin{document}

\title{Morse-Smale surfaced diffeomorphisms with orientable heteroclinic}

\author{A. Morozov}
 \address{National Research University Higher School of Economics, Nizhnij Novgorod, Russia.\\
 ORCID: http://orcid.org/0000-0003-3125-1825}
 \email{andreifrostnn@gmail.com}

\author{O. Pochinka}
 \address{National Research University Higher School of Economics, Nizhnij Novgorod, Russia.\\
 ORCID: http://orcid.org/0000-0002-6587-5305}
 \email{olga-pochinka@yandex.ru}

\begin{abstract} In the present paper we consider preserving orientation Morse-Smale diffeomorphisms on surfaces. Using the methods of factorization and linearizing neighborhoods we prove that such diffeomorphisms have a finite number of orientable heteroclinic orbits.

\end{abstract}

\subjclass{53C17, 22E30, 49J15}

\keywords{Morse-Smale diffeomorphism, factor space, linearizing neighbourhood, heteroclinic intersection}

\maketitle

\section{Introduction} The Morse-Smale dynamics is a pattern of a regular structurally stable (rough) behaviour. The foundations of the theory of roughness of flows, laid down in the classic work by A. Andronov and L. Pontryagin \cite{AnPo}, were developed by the associates of academician A. Andronov in the Gorky school of non-linear oscillations E. Leontovich and A. Mayer \cite{LeMa1}, \cite{LeMa2}. Moreover, A. Meyer introduced the concept of roughness for a discrete dynamical system and obtained a complete topological classification of such systems (Morse-Smale systems) on the circle \cite{Ma}. 

Since then, the theory of topological classification of Morse-Smale dynamic systems has gained wide popularity and has undergone intensive development (see, for example, review \cite{GrGuPoZh}). So the topological classification of structurally stable flows on surfaces is exhaustively described in the work by M. Peixoto \cite{Pe}. On three-dimensional manifolds, the necessary and sufficient conditions of the topological conjugacy of Morse-Smale flows follows from the works by Ya. Umanskiy  \cite{Um}  and A. Prishlyak \cite{Pr}. There are also multidimensional classification results, for example, the classification by S. Pilyugin \cite{Pi} for Morse-Smale flows without heteroclinic intersections on the n-dimensional sphere. For three-dimensional Morse-Smale diffeomorphisms, the completed classification results were obtained only in 2018 by C. Bonatti, V. Grines, O. Pochinka \cite{GrPoBo} after twenty years of work by a team of Russian-French scientists, C. Bonatti, V. Grines, F. Laudenbach, V. Medvedev, E. Pecu, O. Pochinka. The researchers of the Nizhny Novgorod school V. Grines, E. Gurevich, V. Medvedev, O. Pochinka also obtained a complete topological classification of Morse-Smale diffeomorphisms with saddles of co-dimension one on closed n-manifolds \cite{GrGuPoMe}. 

The classification of Morse-Smale cascades on surfaces was obtained in the work by C. Bonatti and R. Langevin \cite{BoLa}, as part of the study of diffeomorphisms with zero-dimensional basic sets. Therefore, it uses a rather heavy apparatus of Markov partitions, due to the need to track the geometry of heteroclinic intersections that make up an infinite number of orbits. Given the obvious evolution of the theory of topological classification of dynamical systems, it would be natural to expect a separate and clear  classification of two-dimensional Morse-Smale casacdes. Such attempts were made for some meaningful subclasses of surface diffeomorphisms, for example, in the works by A. Bezdenezhnykh, V. Grines, T. Mitryakova, O. Pochinka \cite{BeGr}, \cite{GrPoMi}, with the restriction imposed in the form of a finite number of heteroclinic orbits ($beh=1$). 

However, the topology of heteroclinic intersections sometimes itself naturally imposes restrictions on the length of heteroclinic chains, which makes it possible to obtain a classification of such systems in combinatorial terms. In this paper, it is proven that Morse-Smale surfaced diffeomorphism with  orientable heteroclinic intersections cannot have an infinite number of heteroclinic orbits. In the other words we state the following result.

\begin{theorem}\label{thp}  If $f: M^{2}\rightarrow{M^{2}}$ is a preserving orientation Morse-Smale diffeomorphism on surface with orientable heteroclinic then $beh(f)=1$.
\end{theorem}

This result was announced by V. Grines and A. Bezdenezhnyh here \cite{GB} and was proved as the chapter of the candidate's dissertation by A. Bezdenezhnyh \cite{Be}. That proof used a technique of the paper by A. Mayer  \cite{Ma1} about the relation of the number of pairwise disjoint non-trivial recurrent trajectories with the genus of the ambient surface. In this paper we conduct an independent proof of this result using a factorization methods to compatible systems of neighborhoods.

\section{General properties of Morse-Smale diffeomorphisms}
Everywhere below $M^n$ is a smooth  closed connected orientable $n$-dimensional manifold ($n\geq{1}$) and $f:M^n\to M^n$ is a preserving orientation diffeomorphism.

 \begin{definition}[{Morse-Smale diffeomorphism}] Diffeomorphism $f: M^{n}\rightarrow M^{n}$  is called \textit{Morse-Smale diffeomorphism}, if
 \begin{enumerate}
 \item its non-wandering set $\Omega_{f}$ consists of a finite number of hyperbolic orbits;
 \item the invariant manifolds $W^{s}_{p}$, $W^{u}_{q}$ intersect transversally for any non-wandering points $p,q$. 
 \end{enumerate}
 \end{definition}
 
It is well known that Morse-Smale diffeomorphisms have no cycles in the sense of the following definition. 

\begin{definition}[{$k$-cycle}]  A sequence of different periodic orbits  $\mathcal{O}_{0},\mathcal{O}_{1},\dots,\mathcal{O}_{k}$ that satisfies the condition  $$W^u_{\mathcal{O}_{i+1}}\cap W^s_{\mathcal{O}_{i}}\neq\emptyset,i\in\{0,\dots,k\},\mathcal{O}_{k+1}=\mathcal{O}_{0},$$ 
is called \textit{$k$-cycle} ($k\geq 1$).
 \end{definition}

The absence of cycles allows us to introduce the following partial order relation on the set of periodic orbits of the Morse-Smale diffeomorphism (it was first done by S. Smale \cite{Sm}). 
 
\begin{definition}[{Partial order relation, beh}] Let $\mathcal{O}_{i},\mathcal{O}_{j}$ be periodic orbits of the Morse-Smale diffeomorphism $f:M^{n}\rightarrow{M^{n}}$. One says that the orbit $\mathcal{O}_{i},\mathcal{O}_{j}$ are in a relation  $\prec$ ($\mathcal{O}_{i}\prec{\mathcal{O}_{j}}$), if  $${W^{s}_{\mathcal{O}_{i}}}\cap{W^{u}_{\mathcal{O}_{j}}}\neq{\emptyset}.$$
A sequence composed by different periodic orbits $\mathcal{O}_{i}=\mathcal{O}_{i_{0}},\mathcal{O}_{i_{1}},...,\mathcal{O}_{i_{k}}=\mathcal{O}_{j}(k \geq {1})$ such that $\mathcal{O}_{i_{0}}\prec\mathcal{O}_{i_{1}}\prec\ldots\prec\mathcal{O}_{i_{k}}$ is called \textit{chain of length $k$, connecting periodic orbits $\mathcal{O}_{i}$ and $\mathcal{O}_{j}$}. Number equal to the length of the maximum saddle chain of a Morse-Smale diffeomorphism $f:M^n\to M^n$ is denoted $$beh(f).$$ 
\end{definition}

For a subset $P$ of the periodic orbits of $f$ let us set $$beh(\mathcal O_j|P)=\max\limits_{\mathcal O_i\subset P}\{beh(\mathcal O_j|\mathcal O_i)\}.$$

\begin{proposition}[\cite{GrPo}, Theorem 2.1.1] \label{M-Sm-bas} Let $f:M^n\to M^n$ be a Morse-Smale diffeomorphism. Then
\begin{enumerate}
\item  $M^n=\bigcup\limits_{p\in\Omega_f}W^u_p$;
\item  $W^u_p$ is a smooth submanifold of a manifold $M^n$ diffeomorphic to $\mathbb R^{\dim~W^u_p}$ for any periodic point
 $p\in\Omega_f$;
\item $cl(\ell^u_p)\setminus
(\ell^u_p\cup p)=\bigcup\limits_{r\in\Omega_f:\ell^u_p\cap W^s_r\neq\emptyset}W^u_r$ for any unstable separatrix $\ell^u_p$ (a connected component of $W^u_p\setminus p$) of periodic point $p\in\Omega_f$.
\end{enumerate}
\end{proposition}

A similar theorem holds for the stable manifolds of the diffeomorphism $f$.

\section{Surfaced Morse-Smale diffeomorphisms}
In this section we consider a class $MS(M^2)$ of orientation preserving Morse-Smale diffeomorphisms given on a closed orientable surface $M^2$. 
\subsection{Orientable heteroclinic} Following to \cite{Be}, for the diffeomorphisms in this class we will introduce the concept of orientable heteroclinic as the following.

\begin{figure}[h!]
\center{\includegraphics[width=0.8\linewidth]{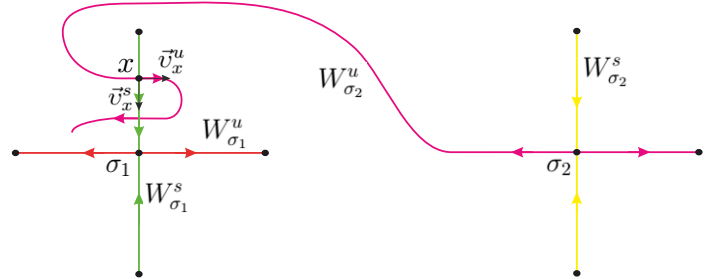}}
\caption{Non-orientable heteroclinic intersection}
\label{ris:or1}
\end{figure}

\begin{figure}[h!]
\center{\includegraphics[width=0.8\linewidth]{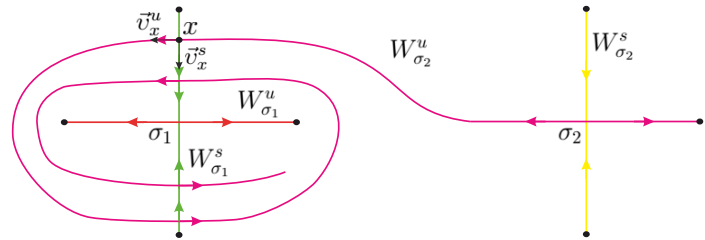}}
\caption{Orientable heteroclinic intersection}
\label{ris:ornot1}
\end{figure}

\begin{definition}[{Orientable heteroclinic}] \label{origp} Let $f\in MS(M^2)$, $\sigma_{i},\sigma_{j}$~--- saddle points of diffeomorphism $f$, such that ${W^{s}_{\sigma_{i}}}\cap{W^{u}_{\sigma_{j}}}\neq\emptyset$. For any heteroclinic point $x\in{W^{s}_{\sigma_{i}}}\cap{W^{u}_{\sigma_{j}}}$ we define an ordered pair of vectors $(\vec{\upsilon}^{u}_{x},\vec{\upsilon}^{s}_{x})$, where:
 \begin{itemize}
 \item  $\vec{\upsilon}^{u}_{x}$ is the tangent vector to the unstable manifold of the point $\sigma_{j}$ at the point $x$ and directed from $x$ to $f^{m^u}(x)$, where $m^u$ is a period of the unstable separatrix containing $x$;
 \item  $\vec{\upsilon}^{s}_{x}$ is the tangent vector to the stable manifold of the point $\sigma_{i}$ at the point $x$ and directed from $x$ to $f^{m^s}(x)$, where $m^s$ is a period of the stable separatrix containing $x$.
\end{itemize}
Heteroclinic intersection of diffeomorphism $f$ is called  \textit{orientable}, if ordered pairs of vectors $(\vec{\upsilon}^{u}_{x},\vec{\upsilon}^{s}_{x})$ determines the same orientation of the ambient surface $M^2$ at every heteroclinic point $x$ of the diffeomorphism $f$. Otherwise, the heteroclinic intersection is called  \textit{non-orientable}  (see Fig. \ref{ris:or1}, Fig. \ref{ris:ornot1}).
\end{definition}

\subsection{Linearizing neighborhoods}

For $\nu\in\{-1,1\}$ denote by $a_{\nu}:\mathbb R^2\to\mathbb R^2$ diffeomorphism given by formula $$a_{\nu}(x_1,x_2)=(\nu\cdot\frac{x_1}{2},\nu\cdot 2x_2).$$ We call $a_{\nu}$ {\it canonical diffeomorphism}. Diffeomorphism $a_{\nu}$ has a single fixed saddle point at the origin $O$ with stable manifold $W^s_O=Ox_1$ and unstable manifold $W^u_O=Ox_{2}$. 

Let  $\mathcal N=\{(x_1,x_2)\in\mathbb{R}^2~:~ |x_1x_2|\leq 1\}$. Notice that the set $\mathcal N$ is invariant with respect to the canonical diffeomorphism
$a_{\nu}$. Define in the neighborhood of $\mathcal N$ a pair of transversal foliations $\mathcal{F}^u,~\mathcal{F}^s$ as follows:   

$$\mathcal{F}^s=\bigcup\limits_{c_{2}\in Ox_2}\{(x_1,x_2)\in \mathcal N~:~x_{2}=c_{2}\},$$ 
$$\mathcal{F}^u=\bigcup\limits_{c_1\in Ox_1}\{(x_1,x_2)\in \mathcal N~:~x_1=c_1\}.$$

Notice that the canonical diffeomorphism $a_{\nu}$ which sends the leaves of foliation $\mathcal{F}^u$ ($\mathcal{F}^s$) to leaves of the same foliation.

Now let $\sigma$ be a saddle point of a diffeomorphism $f\in MS(M^2)$. Suppose $\sigma$ has a period $m_\sigma$ and a type of orientation $\nu_\sigma$, that is $\nu_\sigma=1\,(\nu_\sigma=-1)$ if $f^{m_\sigma}|_{W^u_\sigma}$ preserves (changes) the orientation.
 
\begin{definition}[{Linearizing neighbourhood}] \label{adop} An $f^{m_{\sigma}}$-invariant neighbourhood $N_\sigma$ of saddle point $\sigma\in\Omega_f$ is called by \textit{linearizing} if there is a homeomorphism  ${\mu}_\sigma:N_\sigma\to {\mathcal N}$, conjugate diffeomorphism $f^{m_{\sigma}}\vert_{{N}_\sigma}$ with canonical diffeomorphism $a_{\nu_\sigma}|_{\mathcal N}$.
\end{definition}

\begin{figure}[h!]
\center{\includegraphics[width=0.8\linewidth]{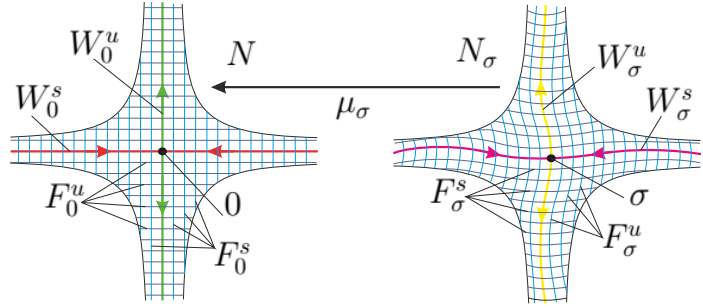}}
\caption{A linearizing neighborhood}
\label{ris:linear_0}
\end{figure}

 Foliations $\mathcal{F}^u,~\mathcal{F}^s$ induced by means of the homeomorphism ${\mu}_\sigma^{-1}$, 
$f^{m_{\sigma}}$-\textit{invariant foliations} ${F}^u_\sigma,~{F}^s_\sigma$ 
on the linearizing neighboughood $N_\sigma$ (see Fig. \ref{ris:linear_0}).

Neighbourhood $N_{\mathcal O_\sigma}=\bigcup\limits_{k=0}^{m_{\sigma}-1}f^k(N_\sigma)$, equipped with a mapping $\mu_{\mathcal O_\sigma}$, composed of homeomorphisms $\mu_\sigma f^{-k}:f^k(N_\sigma)\to \mathcal N,~k=0,\dots,m_{\sigma}-1$, called by \textit{linearizing neighbourhood of orbit $\mathcal O_\sigma$}.

\begin{figure}[h!]
\center{\includegraphics[width=0.8\linewidth]{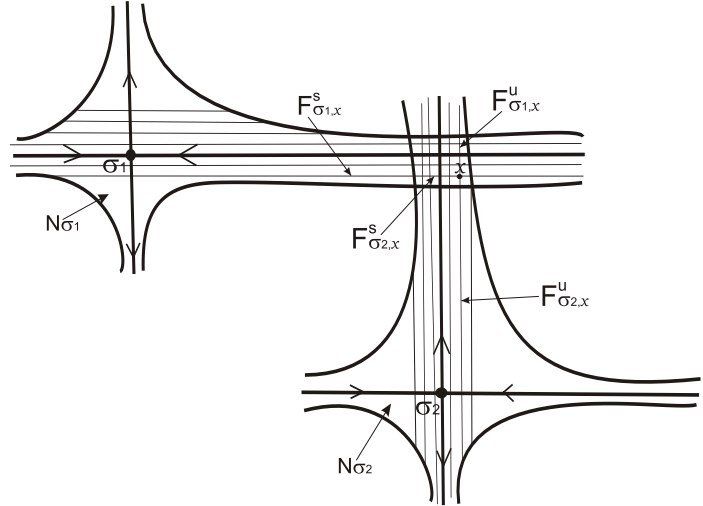}}
\caption{A compatible system of neighborhoods}
\label{ris:neib2}
\end{figure}

\begin{definition}[{The compatible system of neighbourhoods}]
An $f$-invariant collection ${N}_f$ of linearizing neghbourhoods ${N}_\sigma$ of all saddle point $\sigma$ of diffeomorphism $f\in MS(M^2)$ is called \textit{compatible} if the following properties are hold:

\begin{itemize}

 \item if $W^{s}_{{\sigma_1}}\cap W^{u}_{{\sigma_2}}=\emptyset$ and $W^{u}_{{\sigma_1}}\cap W^{s}_{{\sigma_2}}=\emptyset$ for saddle points $\sigma_1,\sigma_2$ then ${{N}}_{{\sigma_1}}\cap{{N}}_{{\sigma_2}}=\emptyset$;

\item if $W^s_{\sigma_1}\cap W^u_{\sigma_2}\neq\emptyset$ for different saddle points $\sigma_1,\sigma_2$ then
$$(F^u_{\sigma_1,x}\cap{N}_{\sigma_{2}})\subset F^u_{\sigma_2,x},\,\,\,(F^s_{\sigma_2,x}\cap N_{\sigma_{1}})\subset F^s_{\sigma_1,x},$$
for  $x\in({N}_{{\sigma_1}}\cap {N}_{{\sigma_2}})$ (see Fig. \ref{ris:neib2}).
\end{itemize}
\end{definition}

\begin{proposition}[\cite{BVPG}, Theorem 1] For every diffeomorphism $f\in{MS(M^{2})}$ there is a compatible system of neighbourhoods. 
\end{proposition}

\subsection{Factorization method} 
Let $f\in MS(M^2)$. Denote by $\Sigma_0,\Sigma,\Sigma_{beh(f)}$ the set of all sinks, saddles, sources of $f$, accordingly. We decompose the set $\Sigma$ on the subsets $\Sigma_1,\dots,\Sigma_{beh(f)-1}$ inductively, as follows: define $\Sigma_i$ as the set of all saddle points of the diffeomorphism $f$, such that $beh(\mathcal O_\sigma|\Sigma_{i-1})=1$ for every orbit $\mathcal O_\sigma,~\sigma\in\Sigma_i$.   

For every  $i\in\{0,\dots,beh(f)-1\}$ define $$\mathcal A_i=\bigcup\limits_{j=0}^{i} W^u_{\Sigma_j},~~\mathcal V_i=\bigcup\limits_{j=0}^{i}W^s_{\Sigma_j}\setminus\mathcal A_i.$$ It follows from \cite{GrPo} that  $\mathcal A_i$ is an attractor of the  diffeomorphism $f$ and $f$ acts discontinuously on the  $\mathcal V_i $. Let $\hat{\mathcal V}_i=\mathcal V_i/f$ and denote by  $$p_{_{i}}:\mathcal V_{i}\to\hat{\mathcal V}_{i}$$ the natural projection.  

Let $\hat V$ be a connected component of $\hat{\mathcal V}_i$ and $V=p_i^{-1}(\hat V)$. Denote by $m_{V}$ the number of connected components in $V$. Set $p_{_V}=p_{_i}|_{V}$. 

\begin{proposition}[\cite{GrPo} Proposition 2.1.5] \label{sos1} The space $\hat V$ is diffeomorphic to two-dimensional torus, natural projection $p_{_V}:V\to\hat V$ is a covering map inducing an epimorphism  $\eta_{_V}:\pi_1(\hat V)\to m_{V}\mathbb Z$ $($here $m_{_V}\mathbb Z$ is  the group of integers multiple to $m_{_V})$ with the following property. Let $[\hat c]\in\pi_1(\hat V)$. Any lift $c$ of the loop $\hat c$ joints a point $x\in V$ with the point $f^n(x)$, where $n\in\mathbb Z$ does not depend on the choice of the lift. Then $\eta_{_V}([\hat c])=n$. 
\end{proposition}

For an unstable separatrix $\gamma^u$ of a saddle point $\sigma$ denote by $m_{\gamma^u}$ its {\it period}, that is, the smallest natural number $\mu$, such that $f^{\mu}(\gamma^u)=\gamma^u$. Let us suppose, that an unstable saddle separatrix  $\gamma^u$ belongs to $V$. Let  $\hat\gamma^u=p_{_V}(\gamma^u)$ and denote by $j_{\hat\gamma^u}:\hat\gamma^u\to\hat V$ the inclusion map. 

\begin{proposition}[\cite{GrPo} proposition 2.1.3] \label{ll} The set  $\hat\gamma^u$ is a circle smoothly embedded in  $\hat V$, such that  $\eta_{_{V}}(j_{\hat\gamma^u*}(\pi_1(\hat\gamma^u)))=m_{\gamma^u}\mathbb Z$. 
\end{proposition}
\begin{figure}[h]
\centerline{\includegraphics[scale=0.2]{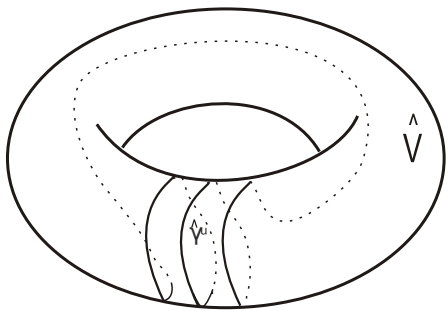}} 
\caption{Projection of an unstable saddle separatrix to $\hat V$}
\label{ris:orbit1}
\end{figure}
On figure \ref{ris:orbit1} the torus $\hat V$ is represented with projection $\hat{\gamma^{u}}$ of separatrix $\gamma^{u}$ such that $\frac{m_{\gamma^{u}}}{m_\omega}=3$.

Let $N_{\gamma^u}$ be a connected component of the set $N_{\sigma}\setminus W^s_{\sigma}$ containing $\gamma^u$ and $\hat N^u_{\gamma^u}=p_{_V}(N_{\gamma^u})$. 

\begin{proposition}[\cite{GrPo} Proposition 2.1.3] \label{l2} The set  $\hat N_{\gamma^u}$ is an annulus smoothly embedded in $\hat V$. 
\end{proposition}

Similar statements can be formulated for the stable saddle separatrix when it belongs to $V$.   

\begin{figure}[h!]
\center{\includegraphics[width=1\linewidth]{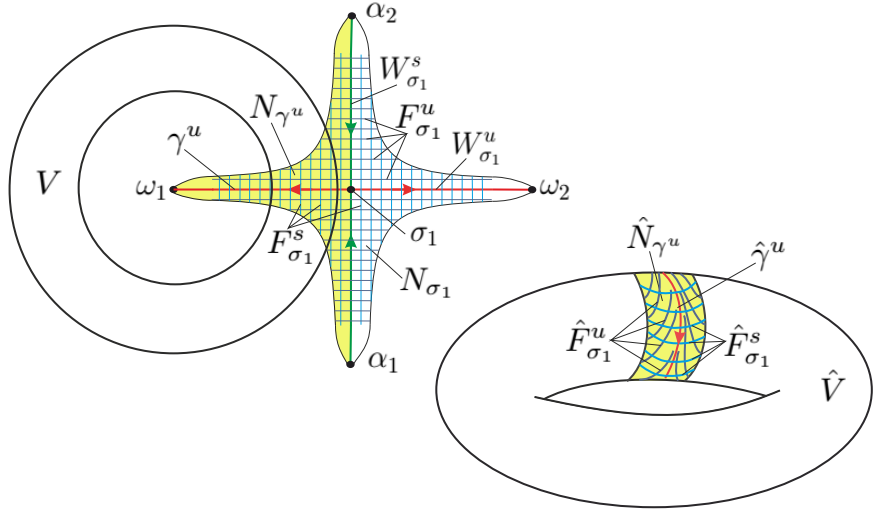}}
\caption{Projection of compatible foliations to $\hat V$}
\label{ris:fs_1sh2}
\end{figure}\par 

On Figure \ref{ris:fs_1sh2} a part of the phase portrait of a diffeomorphism $f\in MS(M^2)$ is represented. Here  $\alpha_{1},\alpha_{2}$ are source points, $\omega_{1},\omega_{2}$ are  sink points, $\sigma_{1}$ is a saddle point and $N_{\sigma_{1}}$ is a linearizing neighbourhood of it. All non-wandering point are supposed to be fixed and $V=W^{s}_{\omega_{1}}\setminus{\omega_{1}}$ --- the punctured basin of the sink $\omega_{1}$. Then on the space orbit $\hat{V}$ we can see annulus $\hat{N}_{\gamma^u} = p_{_V}(N_{\gamma^u})$ as the projection  of a part of $N_{\sigma_{1}}$ containing the unstable separatrix $\gamma^u$. Also we can see projections of both foliations in $N_{\sigma_1}$: $$\hat{F}^{u}_{\sigma_{1}} = p_{_V}(F^{u}_{\sigma_{1}}),\,\,\hat{F}^{s}_{\sigma_{1}} = p_{_V}(F^{s}_{\sigma_{1}}).$$

\section{Proof of the main result (Theorem \ref{thp})} 
Let $f\in MS(M^2)$ and ${N}_f$ be its compatible system of linearizing neghbourhoods ${N}_\sigma$ of all saddle point $\sigma$. In the denotations of previous section for $i=1,\dots,beh(f)-1$ let $$\mathcal W^s_i=\bigcup\limits_{\sigma\in\Sigma_i}W^s_\sigma,\,\,\mathcal W^u_i=\bigcup\limits_{\sigma\in\Sigma_i}W^u_\sigma,$$ $$\mathcal N_i=\bigcup\limits_{\sigma\in\Sigma_i}N_\sigma,\,\,\mathcal F^s_i=\bigcup\limits_{\sigma\in\Sigma_i}F^s_\sigma,\,\,\mathcal F^u_i=\bigcup\limits_{\sigma\in\Sigma_i}F^u_\sigma.$$ Let us prove, that from the conditions of the orientability  heteroclinic intersections of $f$ follows that $beh(f)=1$. 

\begin{figure}[h!]
\center{\includegraphics[width=0.8\linewidth]{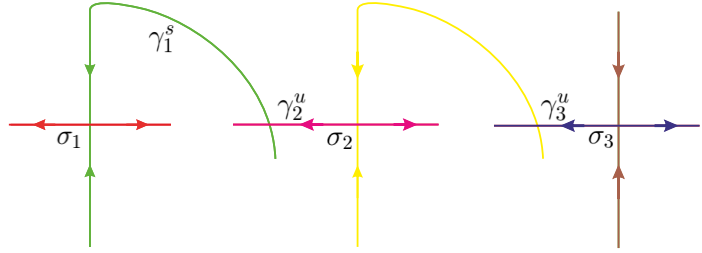}}
\caption{The chain of the periodic points $\sigma_1\prec\sigma_2\prec\sigma_3$ of  the diffeomorphism $f$.}
\label{gamma1+}
\end{figure}

Suppose the opposite: $beh(f)>1$. Then there is a chain of saddle points $\sigma_1\prec\sigma_2\prec\sigma_3$ (see Figure \ref{gamma1+}), such that $\sigma_i\in\Sigma_i$ and (see Fig. \ref{gamma1+}) $$W^u_{\sigma_{i+1}}\cap W^s_{\sigma_i}\neq\emptyset,i=1,2.$$ 

Then these saddle points have separatrix $\gamma^s_1,\gamma^u_2,\gamma^u_3$ such that $\gamma^s_1\cap\gamma^u_2\neq\emptyset$ and $\gamma^u_3\cap W^s_{\sigma_2}\neq\emptyset$. Let $\hat \gamma^s_{1}=p_1(\gamma^s_{1}),\,\hat\gamma^u_{2}=p_1(\gamma^u_{2})$. By the construction, all points of the separatrix $\gamma^s_{1},\gamma^u_{2}$ belong to space $\mathcal V_1$. Due to Proposition \ref{ll}, each set $\hat {\gamma}^s_{1}, \hat\gamma^u_{2}$ is a circle on a disjoint union of the tori $\hat{\mathcal V}_1$. Since the intersection $\hat\gamma^s_{1}\cap\hat\gamma^u_{2}$ is not empty then these circles belong to the same torus (denote it by $\hat V$). Also we will see annuli $\hat N_1=p_1(N_{\gamma^s_1})$ and $\hat N^1_2=p_1(N_{\gamma^u_2})$ on the torus $\hat V$.

Every point of the separatrix $\gamma^u_3$, except the finite number of heteroclinic orbits $\gamma^u_3\cap \mathcal W^s_{2}$, belongs to the space $\mathcal V_1$. Then the set $\hat\gamma^u_3=p_1(\gamma^u_3)$ consists of a finite number of non-compact arcs. As $\gamma^u_3\cap W^s_{\sigma_2}\neq\emptyset$ then, due to Proposition \ref{M-Sm-bas}, $\gamma^u_3\cap N_{\gamma^u_2}\neq\emptyset$. Thus there is a connected component $\hat\ell$ of the set $\hat\gamma^u_3$ such that $\hat\ell\cap\hat N^1_2\neq \emptyset$.  

Let $\hat\Gamma^u_2=\hat V\cap p_1(\mathcal W^u_2)$ and $\hat{\mathcal N}_2=p_1(\mathcal N_{2})$. Due to Proposition \ref{ll} the set $\hat\Gamma^u_2$ is a disjoint union of homotopically non-trivial circles on the torus $\hat V$ and $\hat\gamma^u_2$ is one of them. From  condition of the orientability of  heteroclinic  intersections, it follows that the circle $\hat\gamma^s_1$ intersects each of the circles of the set $\hat\Gamma^u_2$ and the index of the intersection at each point of  such the intersection is the same. 

\begin{figure}[h!]
\center{\includegraphics[width=0.8\linewidth]{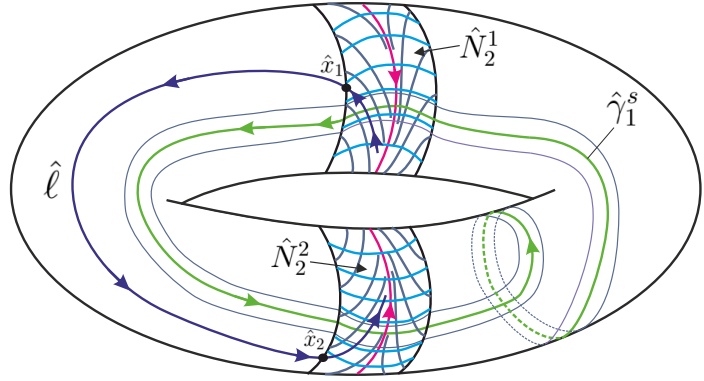}}
\caption{Orbit space $\hat V$}
\label{finalfroof}
\end{figure}

As the curve $\hat\ell$ belongs $\hat V$ then there is annuli $\hat N^2_2$ (possible $\hat N^2_2=\hat N^1_2$) which is a connected components of the set $\hat{N}_2$ and such that $(\hat\ell\cap\hat{\mathcal N}_2)\subset (\hat N^1_2\cup\hat N^2_2)$. From the compatible system of neighborhoods follows, that each connected component of the intersection $\hat\ell\cap(\hat N^1_2\cup\hat N^2_2)$ is a leaf of the foliation $\hat{\mathcal F}^u_2=p_1(\mathcal F^u_2)$ (see Fig. \ref{finalfroof}). 

Let $\hat x_1=\hat{\ell}\cap{\partial{\hat{N}^{1}_{2}}}$ and $\hat x_2=\hat{\ell}\cap{\partial{\hat{N}^{2}_{2}}}$. Due to the orientability of the heteroclinic intersection the points $\hat x_1,\hat x_2$ belong to different connected components of the set $\partial (\hat N^1_2\cup \hat N^2_2)$. Moreover, according to compatibility of system of linearizing neighborhoods,  the arc $\hat\ell$ is a union of three arcs with the pairwise disjoint interiors $$\hat\ell= \hat F^u_{1,{\hat{\ell}}}\cup\hat\ell_{\hat x_1,\hat x_2}\cup  \hat F^u_{2,{\hat{\ell}}},$$ where $ \hat F^u_{1,{\hat{\ell}}},\hat  F^u_{2,{\hat{\ell}}}$ are different leaves of foliation  $\hat{\mathcal F}^u_2$. 

\begin{figure}[h!]
\begin{minipage}[h!]{0.45\linewidth}
\center{\includegraphics[width=1\linewidth]{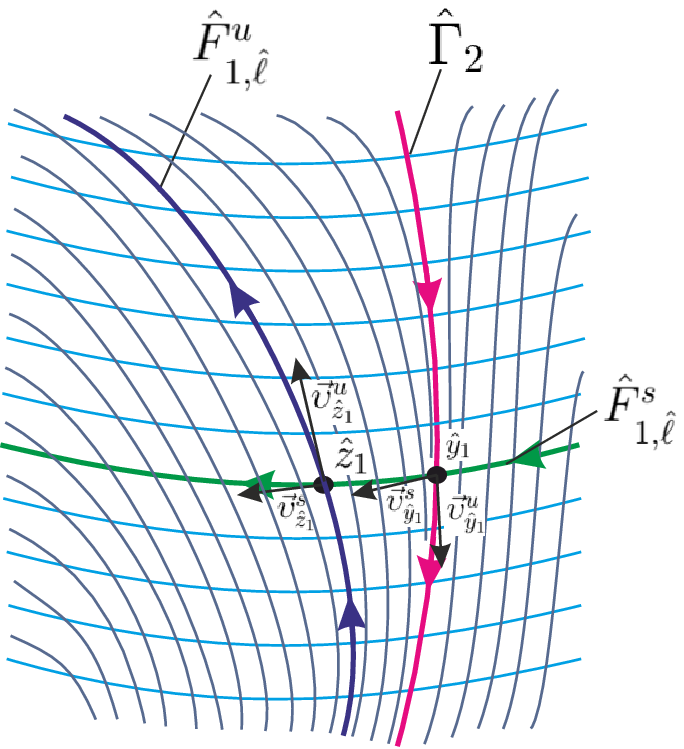} \\ a)}
\end{minipage}
\hfill
\begin{minipage}[h!]{0.45\linewidth}
\center{\includegraphics[width=1\linewidth]{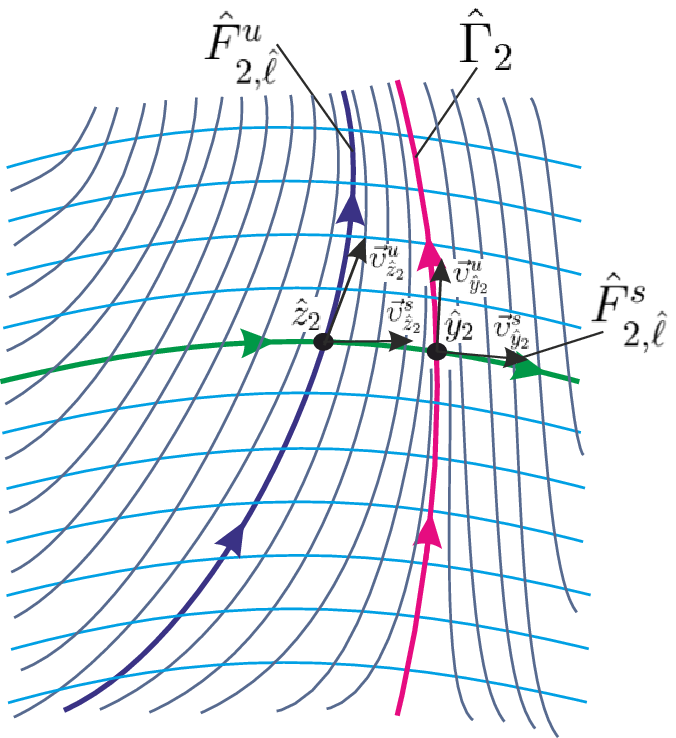} \\ b)}
\end{minipage}
\caption{Intersection of separatrices in the space orbit $\hat{V}$}
\label{n12}
\end{figure}

It follows from the definition of the compatible system of linearizing neighborhoods that every connected component of the intersection $\hat{\gamma}^s_1\cap(\hat{N}^{1}_{2}\cup\hat{N}^{2}_{2})$ is a leaf of the foliation $\hat{\mathcal F}^s_2=p_1(\mathcal F^s_2)$. Choose one connected component $\hat F^s_{1,{\hat{\ell}}}\,(\hat  F^s_{2,{\hat{\ell}}})$ in the intersection $\hat{\gamma}^s_1\cap\hat{N}^{1}_{2}\,(\hat{\gamma}^s_1\cap\hat{N}^{2}_{2})$, possible $\hat F^s_{1,{\hat{\ell}}}=\hat  F^s_{2,{\hat{\ell}}}$ if $\hat N^1_2=\hat N^2_2$. Let $\hat y_1=\hat F^s_{1,{\hat{\ell}}}\cap \hat\Gamma_2$ and $\hat y_2=\hat F^s_{2,{\hat{\ell}}}\cap \hat\Gamma_2$. Also choose one pair of the points  $\hat z_1,\hat z_2$ from each of the countable set $\hat F^s_{1,{\hat{\ell}}}\cap \hat F^u_{1,{\hat{\ell}}}$, $\hat F^s_{2,{\hat{\ell}}}\cap \hat F^u_{2,{\hat{\ell}}}$ (see Fig. \ref{n12}). 

According to the orientability of the heteroclinic intersection, the pairs of the vectors $(\vec{\upsilon}^{u}_{\hat y_{1}},\vec{\upsilon}^{s}_{\hat y_{1}})$ and $(\vec{\upsilon}^{u}_{\hat y_{2}},\vec{\upsilon}^{s}_{\hat y_{2}})$ define the same orientation of the torus $\hat V$. Also, the pairs of the vectors  $(\vec{\upsilon}^{u}_{\hat y_{1}},\vec{\upsilon}^{s}_{\hat y_{1}})$, $(\vec{\upsilon}^{u}_{\hat z_{1}},\vec{\upsilon}^{s}_{\hat z_{1}})$ and $(\vec{\upsilon}^{u}_{\hat y_{2}},\vec{\upsilon}^{s}_{\hat y_{2}})$, $(\vec{\upsilon}^{u}_{\hat z_{2}},\vec{\upsilon}^{s}_{\hat z_{2}})$ also need to be consistently oriented. However, this contradicts with any orientation on the curve $\hat\ell$.

\newpage

\subsection*{Acknowledgements.} The proof of Theorem \ref{thp} is supported by RSF (Grant no. 17-11-01041), the proof of auxiliary results is supported by Basic Research Program at the National Research University Higher School of Economics (HSE) in 2019.

\end{document}